\documentclass[12pt,reqno]{amsart}   	
\usepackage[letterpaper, margin=1in]{geometry}
\usepackage{amssymb}
\usepackage{amsmath}
\usepackage{amsfonts}
\usepackage[utf8]{inputenc}
\usepackage{amsthm}
\usepackage{tikz-cd}
\usepackage{array}

\newtheorem{thm}{Theorem}[section]
\newtheorem{lemma}[thm]{Lemma}
\newtheorem{df}[thm]{Definition}
\newtheorem{prop}[thm]{Proposition}
\newtheorem{ex}[thm]{Example}
\newtheorem{cor}[thm]{Corollary}

\theoremstyle{obs}

\theoremstyle{prob}
\newtheorem{prob}[thm]{Problem}
\theoremstyle{conj}
\newtheorem{conj}[thm]{Conjecture}

\newtheoremstyle{remark}
    {\dimexpr\topsep/2\relax} 
    {\dimexpr\topsep/2\relax} 
    {}          
    {}          
    {\bfseries} 
    {.}         
    {.5em}      
    {}          

\theoremstyle{remark}
\newtheorem{remark}[thm]{Remark}
\newcommand{\RR}{\mathcal{R}}
\newcommand{\Z}{\mathbb{Z}}

\newcommand{\R}{\mathbb{R}}

\newcommand{\OO}{\mathcal{O}}
\newcommand{\C}{\mathbb{C}}

\title{Sums of squares of regular functions on rational surfaces}
\author{Tomasz  Kowalczyk}
\date{}

\begin{document}

\maketitle

\begin{abstract}
We study the sums of squares on cylinders of the form $X \times \mathbb{A}_K$ for a (weakly) factorial curve $C$. We prove the equality of the Pythagoras numbers of the ring of regular functions on the cylinder with that of the field of rational functions. We then apply these results to the case of (uniformly) rational varieties. We show that if $X$ is a nonsingular rational algebraic surface over the reals, then the Pythagoras number of the ring of regular functions on $X$ is bounded above by 12.

\end{abstract}

\section*{Introduction}

 Let $f$ be a real valued function defined on a space $X$. Assume that $f$ is a sum of squares of functions on $X$. What is the minimal number of squares needed to represent $f$? In this paper we tackle this problem for regular functions on an algebraic surface defined over a formally real field $K$. This problem is closely related to that of computing the Pythagoras number of a ring.

 \begin{df}
For a commutative ring $R$ with identity its Pythagoras number, denoted by $p(R)$, is the smallest positive integer $g$ such that any element which is a sum of squares can be written as a sum of at most $g$ squares. If such number does not exist, we put $p(R)=\infty$.
\end{df}

The problem of computing the Pythagoras number is a generalization of the Waring problem for positive integers. Hence, computing Pythagoras numbers is essentially the Waring problem for squares for more general rings.

For a commutative ring $R$ with identity, determining its Pythagoras number is a rather difficult task. As we will be considering real rings (i.e. rings in which $-1$ is not a sum of squares), let us discuss what is known in this case. Any 1--dimensional affine $\R$--algebra has finite Pythagoras number (see \cite{cldr1982}). On the other hand, any affine $K$--algebra over a formally real field $K$ of dimension at least 3 has infinite Pythagoras number. Similarly, any regular local ring $R$ of dimension at least 3 with formally real residue field has $p(R)=\infty$ (see \cite{cldr1982} for both results). Let us also mention a very recent preprint \cite{osz2025} which deals with similar problems.

The case of dimension 2 seems to be the most mysterious. It is known that $p(\R[[x,y]])=2$ and $p(\Z[[x]])=5$, but $p(\Z[x])=p(\R[x,y])=\infty$ (see \cite{cldr1982}).  Recently, there has been some development in this area. Pythagoras number of rings of the type $\R[x,y,\sqrt{f(x,y)}]$, where $f(x,y)$ is a polynomial satisfying some mild conditions, has been computed in \cite{bk2024}. See also \cite{bk2023} for the $k$-regulous case, and \cite{vill2021} for some other problems related to sums of squares. For a 2--dimensional regular local ring $R$ it was shown in \cite{scheiderer2001} that finiteness of the Pythagoras number of the field of fractions of $R$ implies finiteness of the Pythagoras number of $R$ and explicit bound was given. We are able to strengthen this result in some cases, see Corollary \ref{local}. See also \cite{bz2024} for the case of the field of fractions of formal power series over $\R$.

Computing the Pythagoras number of a number field is a rather simple task \cite[Chapter 7, Example 1.4]{pfister1995}. On the other hand, computing the Pythagoras number of the ring of integers or orders is much more involved, see \cite{dgmy2025, dmr2021, kc2018, krasensky2022, krs2022, ky2023, tinkova2023} for some recent developments in this area. It is also worth to mention very recent paper \cite{pop2023}, which deals with the Pythagoras number of the field of rational functions of a curve over a number field. This result has already been improved in \cite{benoist2025}.

In the literature one can also find treatment of higher Pythagoras numbers (sometimes called Waring numbers) (cf. \cite{becker1982, grimm2015, km2024, kv2023}). Computing Pythagoras number is difficult, and computing higher Pythagoras numbers is much harder, as the theory of higher degree forms is much less developed.

Our main result relates the length of a regular function to the length of this function considered as a rational function. Surprisingly, the Pythagoras number was not known for a ring of regular functions of any surface, not even the real plane. Let us discuss what is known. If $R$ is a "ring of regular functions" (RFR) in the sense of \cite{mahe1990}, then any totally positive unit can be written as a sum of squares, and explicit upper bound was given in terms of transcendence degree. For an algebraic surface over a formally real field $K$, the usual ring of regular functions is an RFR, hence the result of Mah\'e applies. This is however far from computing the Pythagoras number, as there is no control on the lengths of non--units. In a similar spirit one defines the unital Pythagoras number i.e. one considers only units which are sums of squares and not general elements which are sums of squares. It is worth mentioning the paper \cite[Theorem 5.17]{bp1996} where this notion was studied. The most interesting part is the relation with higher unital Pythagoras numbers. Such relation is known only in  very specific cases for usual Pythagoras numbers \cite[Conjecture 5.2]{kv2023}.

Let us describe the content of this paper. 
In Section 1 we recall necessary facts from algebraic geometry. Section 2 contains discussion of (weakly) factorial varieties. We also study sums of squares in rings of regular functions on curves. In the next section we study sums of squares on cylinders i.e. products of a curve and a line. Section 4 is devoted to the study of (uniformly) rational surfaces. Among other results, we show that if $X$ is a smooth nonsingular rational surface over $\R$ then $p(\OO(X)) \leq 12$. We finish the paper with discussion of possible further generalizations of theorem of Cassels, and list some conjectures and open problems which may stimulate further research.



\section{Preliminaries}

Let us begin with some preliminary definitions. We will start with the notion of  length, which will be useful in the study of sums of squares. Let $R$ be a commutative ring with identity.
\begin{df}
 We define the length of an element $a \in R$ to be the smallest positive integer $g$ such that $a$ can be written as a sum of at most $g$ squares. If $a$ is not a sum of squares, then its length is infinite.
\end{df}
\begin{ex}
Let $R=\mathbb{Z}$ and $a=7$. In this case it is known that the length of $7$ is $4$.
\end{ex}

From the very definition of length we obtain the following fact

\begin{lemma}\label{lemma}
Let $R$ be an integral domain and $a \in R$.
Then the length of $a$ in $R$ is greater than or equal to the length of $a$ in the field of fractions of $R$.  \qed
\end{lemma}
This in particular implies that $p(R) \geq p(K)$ where $K$ is the fraction field of $R$.

\begin{df}
We say that a field $K$ is formally real if $-1$ is not a sum of squares of elements of $K$.

\end{df}

In particular, any formally real field $K$ contains $\mathbb{Q}$ and so it is of characteristic $0$.

Let $K$ be a formally real field. We will recall all the definitions regarding algebraic geometry and rings of regular functions. In this paper we will avoid the scheme theoretic point of view on algebraic varieties.
For the treatment of algebraic sets over a real closed field see \cite[Chapter 3]{bcr1998}. Let $X \subset K^n$ be a Zariski closed algebraic set.
\begin{df}
We define the coordinate ring of $X$ as the quotient ring 
$$K[X]:=K[x_1,x_2,\dots, x_n]/\mathcal{I}(X)$$
 where $\mathcal{I}(X)$ is the ideal of all polynomials vanishing on $X$.

If $X$ is an irreducible algebraic set, then we denote field of fractions of $K[X]$ by $K(X)$.
\end{df}

Take any $x_0\in X$ and denote by $\mathfrak{m}_{x_0}$ the maximal ideal of all polynomial functions vanishing on $x_0$.

\begin{df}
We say that $x_0 \in X$ is a nonsingular point if the local ring $K[X]_{\mathfrak{m}_{x_0}}$ is a regular local ring.  
We say that an algebraic set $X\subset K^n$ is nonsingular if every point point is nonsingular.

\end{df}

As we will be interested only in $K$-rational points, we may allow singularities to occur in the algebraic closure of $K$. If $K$ is not algebraically closed field, then the coordinate ring lack many key properties. Usually, in real algebraic geometry, a much more natural ring is the ring of regular functions -- the crucial object of our study.

\begin{df}
Let $X \subset K^n$ be an algebraic subset. We define the ring of regular functions on $X$ as 
$$\OO(X)= \left \{ \,  \frac{f}{g} \,  | \,  f,g,\in K[X], \, g^{-1}(0)=\varnothing  \right\}.$$

\end{df}

Let us note that similarly as in the scheme theory one can consider sheaves of modules over the sheaf of regular functions. Our geometric approach differs substantially from the case of an algebraically closed field (or scheme theoretic point of view). For example, for real algebraic sets with the sheaf of regular functions, Cartan's Theorems $A$ and $B$ do not hold, however they do hold after a suitable multiblowup \cite{kowalczyk2018, kowalczyk2019}.

In algebraic geometry one defines the notion of a regular function locally. However, since we deal only with formally real fields, by \cite[Proposition 3.2.3]{bcr1998} regular functions on $X$ are precisely those that admit a global denominator. In particular, $\OO(X)$ is a localization of $K[X]$ at a set consisting of polynomial functions vanishing nowhere on $X$.



\begin{remark}
By an algebraic variety we will mean an irreducible algebraic set.
\end{remark}

\begin{remark}
In this paper we will consider product surfaces, hence to avoid awkward formulas such as $K[C\times K]$, we will denote the affine line over $K$ as $\mathbb{A}_K$.
 \end{remark}

Let $X \subset K^n$ be an irreducible algebraic set.

\begin{df}
We say that $X$ is a factorial variety if the ring $K[X]$ is a unique factorization ring. 
We say that $X$ is a weakly factorial variety if the ring $\OO(X)$ is a unique factorization ring.
\end{df}

If $X$ is a factorial variety then it is obviously weakly factorial. Also, by a classical result in algebra, if $X$ is factorial, the so is cylinder $X\times \mathbb{A}_K$ over $X$. We will now discuss properties of weakly factorial varieties.

\begin{lemma}
Let $X$ be a weakly factorial variety. Then the cylinder $X\times \mathbb{A}_K$ is a weakly factorial variety.

\begin{proof}
By the assumption $\OO(X)$ is a UFD, hence so is $\OO(X)[y]$ by \cite[Theorem 6.14]{hungerford1980}. Of course $K[X\times \mathbb{A}_K] \subset \OO(X)[y]$. Let $S$ be the multiplicative set of all polynomial functions vanishing nowhere on $X\times \mathbb{A}_K$. Straightforward calculation shows that $S^{-1}K[X\times \mathbb{A}_K]  = \OO(X\times \mathbb{A}_K)=S^{-1}\OO(X)[y]$. As a consequence, $\OO(X\times \mathbb{A}_K)$ is a UFD as a localization of a UFD.
\end{proof}
\end{lemma}
\begin{ex}
It is known that 1-dimensional unit sphere $\mathbb{S}^1$ is not a factorial variety, while 2-dimensional unit sphere $\mathbb{S}^2$ is, both considered as algebraic sets over $\R$ (see \cite[Theorem 5]{swan1962}). The hyperbola $X=\{xy=1 \}\subset K^2$ is also factorial, since its ring of polynomial functions $K[X]\cong K[x,x^{-1}]$ is a localization of $K[x]$.
\end{ex}

\section{Weakly factorial varieties over $\R$.}
If we are working over the field of real numbers, we are able to give equivalent conditions for weak factoriality. In particular, we have 

\begin{thm}\cite[Proposition 12.4.14]{bcr1998} Let $X$ be a nonsingular compact algebraic variety over $\R$ of dimension $d$. Then the following conditions are equivalent:
\begin{itemize}
\item $X$ is weakly factorial,  
\item $H_{d-1}^{alg}(X, \Z/2)=0$.
\end{itemize}

\end{thm}
Here, $H_{d-1}^{alg}(X, \Z/2)$ is the group of $d-1$ homology classes generated by classes represented by compact algebraic subsets. The above Theorem holds, for example if $H_{d-1}(X, \Z/2)=0$. As pointed out by Shiota, the factoriality of $X$ is not a topological invariant as the following example shows.
\begin{ex}\cite[Exemple]{shiota1981}
Consider two algebraic sets, $X= \{ x^2+y^2+z^2=1\}$ and $Y=\{x^4+y^4+z^2=1 \}$. $X$ and $Y$ are homeomorphic, both are weakly factorial and $X$ is factorial while $Y$ is not. To see this, note that in $\R[Y]$ we have $x^4+y^4=1-z^2=(1-z)(1+z)=(x^2+\sqrt{2}xy+y^2)(x^2-\sqrt{2}xy+y^2)$ (cf. \cite[Section 7]{bks1982}).

\end{ex}
Non compact case was also studied, however much less is known (cf. \cite[Theorem 5, Remark 6]{bks1982}). In particular, if $X$ is a nonsingular nonorientable algebraic set, then $X$ is not weakly factorial.

In this paper, weakly factorial curves will be of great importance. Let us now provide a topological classification of those.
\begin{lemma}\cite[Theorem 2.3]{gw1980}
Let $X$ be a nonsingular algebraic curve over $\R$. Assume that $X$ has a compact connected component. Then $X$ is not a weakly factorial variety, hence it is not factorial.
\end{lemma}
In its original statement, the above Theorem shows that the maximal ideal corresponding to a point in a compact connected component cannot be principal. This is shown for the coordinate ring $\R[X]$ as well as $S^{-1}\R[X]$ where $S$ is the multiplicative set of all polynomial functions vanishing nowhere on the given connected component. In particular, this ring is a localization of the ring of regular functions $\OO(X)$, hence $X$ is not weakly factorial.
\begin{thm}
Let $X$ be a nonsingular algebraic curve over $\R$ with no compact connected component. Then $X$ is weakly factorial.

\begin{proof}
First of all, let us note that a nonsingular real algebraic curve is homeomorphic to a disjoint union of circles and intervals. Since $X$ has no compact connected component, $X$ is homeomorphic to a disjoint sum of finitely many open interval, hence $H^1(X, \mathbb{Z}_2)=0$.
Secondly, let us embed $X \subset \mathbb{P}^n(\R)$ for some $n$. Denote by $X'$ the possible desingularization of the Zariski closure of the image of $X$ in $\mathbb{P}^n(\R)$. In this case $X'\setminus X$ consists of finitely many points. Moreover every class in the homology group of degree $0$, i.e. $H_0(X'\setminus X, \Z_2)$ is represented by a compact algebraic subset. The result now follows readily from the theorem below.

\end{proof}

\end{thm}

\begin{thm}\cite[Th\'eor\`eme 2]{shiota1981}
Let $X$ be an irreducible, non compact algebraic set of dimension $n$ with $H^1(X, \Z_2)=0$. If every element of $H_{n-1}(X'\setminus X, \Z_2)$ is realized by a compact algebraic set then $X$ is weakly factorial.
\end{thm}

The above results can be combined into the following.
\begin{thm}\label{class}
Let $X$ be a nonsingular algebraic curve over $\R$. Then $X$ is weakly factorial if and only if $X$ has no compact connected component.

\end{thm}
We see that for a nonsingular irreducible curve $C$ over $\R$, weak factoriality is a local property in the sense that any point $x \in C$ has a Zariski open neighborhood $U$ such that $\OO(U)$ is a weakly factorial curve.

\begin{remark}
The results of this Section essentially use the properties of the real numbers. In would be interesting in its own right to find a generalization of the above to arbitrarily fields, cf. Section \ref{sec:open problems}.
\end{remark}

Let us a finish this section with a discussion of sums of squares of regular functions on curves. Let $K$ be a formally real field.

\begin{prop}\label{p_2O=p_2K}
Let $C \subset K^n$ be an irreducible weakly factorial curve. Then $p(\OO(C))=p(K(C))$.

\begin{proof}
Since $K(C)$ is the field of fractions of $\OO(C)$ it is enough to show that the length of a regular function $f$ is the same in both aforementioned rings.  Let $f \in \OO(C)$ be a sum of squares. By the weak factoriality, we may factor $f$ as $f=g_1^2g$, where $g_1,g \in \OO(C)$ and $g$ is a unit. Write now $g$ as a sum of squares in $K(C)$ i.e. $g=\frac{\sum_{i=1}^Nf_i^2}{h^2}$ with $h,f_i \in K[C]$ for $i=1,2, \dots ,N$. By the weak factoriality, any possible zero of $h$ could be factored out from all of the $f_i$'s. We then obtain a representation of $g$ with a possibly smaller number of summands, all of which are regular functions. Hence the length of $f$ in $\OO(C)$ and $K(C)$ are equal.

\end{proof}

\end{prop}

We may strengthen the above result for $K=\R$.

\begin{thm}\label{p_2(OO(C))}
If $C$ be a nonsingular irreducible algebraic curve over $\R$ then $p(\OO(C)) \leq 4$.

\begin{proof}
The result of Pfister \cite[Corollary 3.4]{pfister1995} implies that $p(\R(C)) \leq 2.$ By the Theorem \ref{class} we may take a covering of $C$ by two open subsets $U_1,U_2 \subset C$ such that $\OO(U_i)$ is a UFD. Let $f \in \OO(C)$ be a sum of squares. Choose functions $h_1,h_2 \in \OO(C)$ such that $\mathcal{Z}(h_i)=C\setminus U_i$.
By the previous proposition, we have that $f_i=f|_{U_i}$ is expressible as a sum of at most two squares in $\OO(U_i)$, i.e. $f_i= \frac{f_{1i}^2+f_{i2}^2}{g_i^2}= \frac{(f_{1i}^2+f_{i2}^2)h_i^2}{g_i^2h_i^2}$, where the zero set of $g_i$ is contained in $C \setminus U_i$ and $f_{ij}\in \OO(U_i)$, $i,j=1,2$. Having those two presentations we may add numerators and denominators together to obtain a global presentation of $f$ i.e.
$$f=\frac{(f_{11}^2+f_{12}^2)h_1^2+(f_{21}^2+f_{22}^2)h_2^2}{g_1^2h_1^2+g_2^2h_2^2}=\frac{\left((f_{11}^2+f_{12}^2)h_1^2+(f_{21}^2+f_{22}^2)h_2^2 \right)(g_1^2h_1^2+g_2^2h_2^2)}{(g_1^2h_1^2+g_2^2h_2^2)^2}.$$ By construction, the denominator does not vanish on $C$. Applying the 2-square (or 4-square) identity in the numerator finishes the proof.
\end{proof}
\end{thm}
In the above proof we have shown the existence of a presentation as a quotient of regular functions, however, we have defined regular functions as quotients of polynomial functions. Since the denominators of global regular functions cannot vanish, one can always simplify such presentations to the desired form if necessary.
\begin{remark}
The above results can be also shown for higher even Pythagoras numbers (cf. \cite{kv2023}). The proofs are essentially the same, hence we will omit them.
\end{remark}

\section{Sums of squares on cylinders.}
As we have discussed the case of curves, the main aim of this section is to generalize Proposition \ref{p_2O=p_2K} to a possibly large class of algebraic surfaces.

Apart from algebraic geometry, we will use elements of the theory of quadratic forms (see \cite{pfister1995}). Let us now recall the Theorem of Cassels \cite{cassels1964}, where only sums of squares were considered instead of a general quadratic form. The statement we present below can be found  in \cite[Chapter 1, Theorem 2.2]{pfister1995}.

\begin{thm}
Let $K$ be a field of characteristic different from 2. Let $\varphi$ be a quadratic forms with coefficients in $K$ and $f \in K[x]$ be a polynomial. If $f$ is represented by $\varphi$ over $K(x)$ then it is represented by $\varphi$ already over $K[x]$.
\end{thm}


 



We are now ready to prove the main result of this Section.
\begin{thm}\label{main}
Let $K$ be a formally real field. Let $C$ be an irreducible nonsingular weakly factorial curve over $K$, and $X=C\times \mathbb{A}_K$ be a product surface. Then $p(\OO(X))=p(K(X))$.

\begin{proof}
Similarly as before, since $K(X)$ is the field of fractions of $\OO(X)$, it is enough to show that the lengths of a regular functions $f$ in $K(X)$ and $\OO(X)$ are equal. By the weak factoriality of $X$, we may assume that the zero set of $f$ is finite and $f$ does not vanish on any curve. Denote the length of $f$ in $K(X)$ by $N$.


Of course, $K[C][y] \subset K(C)[y]$ and let $\varphi = \sum_{i=1}^N x_i^2$ be a quadratic form. By the very definition of $N$, $\varphi$ represents $f$ over $K(X)$, hence, by the Theorem of Cassels also over $K(C)[y]$. Let $f=\sum_{i=1}^N \frac{f_i^2}{g^2}$ be a presentation such that $f_i \in K[X]$ and $g \in K[C]$. If $g$ does not vanish on $C$ then we are done. Any zero $x_0 \in C$ of $g$ corresponds to a curve of the form $\{x_0 \} \times \mathbb{A}_K$ in $X$ and those are given by a single equation in $\OO(X)$ by the weak factoriality of $C$ and $X$. Hence, we can factor them from $g$ and all $f_i$ for $i=1,2,\dots, N$. However, since the factorization takes places in $\OO(X)$, it might happen that after this procedure the denominator $g$ might be replaced by a nowhere vanishing function on the whole $X$, not necessarily in $K[C]$, but in any case, $g$ does not vanish and so we obtain representation of $f$ as a sum of squares of regular functions of length $N$.



\end{proof}
\end{thm}

\begin{remark}
If we further assume that the curve $C$ is factorial, then in the above proof we may assume that the denominator $g$ does not depend on $y$. This follows from the fact, that we are able to cancel out factors in the ring $K[X]$ rather than in $\OO(X)$.

\end{remark}

\begin{cor}\label{cor:p(O(X))=p(O(U))}
With the assumptions as in Theorem \ref{main}. Let  $U \subset X$ be a Zariski open subset of $X$. Then $p(\OO(X))=p(\OO(U))$.
\end{cor}
\begin{proof}
If $K$ is not an algebraically closed field then any Zariski open subset is isomorphic to an algebraic set. In this case $\OO(U)$ is isomorphic with a localization of $\OO(X)$ and so we can derive
$p(\OO(X))\geq p(\OO(U)) \geq p(K(X))$. By the previous theorem the proof is done.
\end{proof}

Two above results immediately imply the following
\begin{cor}
Let $U \subset \R^2$ be a Zariski open subset. Then $p(\OO(U))=p(\OO(\R^2))=4$.
\end{cor}

Similarly as for curves, we may strengthen Theorem \ref{main} for cylinders over $\R$.

\begin{thm}
If $C$ is a nonsingular irreducible algebraic curve over $\R$, then ${p(\OO(C\times \R)) \leq 8}$.

\begin{proof}
It is enough to combine Theorem \ref{main} and the proof of Theorem \ref{p_2(OO(C))} with the result of Pfister  \cite[Corollary 3.4]{pfister1995} which implies that $p(\R(C\times \R)) \leq 4.$
\end{proof}
\end{thm}

Let $R$ be a regular local ring of dimension 2. Denote by $F$ the field of fractions of $R$. Problem 5 in \cite{cldr1982} asked about a relation between $p(R)$ and $p(F)$. Scheiderer \cite{scheiderer2001} proved that $p(R)\leq 4p(F)-4$. In particular, $p(R)$ is finite if and only if $p(F)$ is finite. We are able to strengthen this result, at least for some regular local rings of geometric origin.

\begin{cor}\label{local}
With the assumptions of Theorem \ref{main},  let $x_0 \in X $ be any nonsingular point. Let $\OO(X)_{\mathfrak{m}_{x_0}}$ be the regular ring of the point $x_0$. Then $p(\OO(X)_{\mathfrak{m}_{x_0}})=p(K(X))$.

\begin{proof}
Since $\OO(X)_{\mathfrak{m}_{x_0}}$ is a localization of $\OO(X)$, the result follows readily from Theorem \ref{main}.
\end{proof}

\end{cor}

If $X$ is a surface defined over a real closed field $K$, then it follows from the Pfister's theory of multiplicative forms \cite[Corollary 3.4]{pfister1995} that $p(K(X))\leq 4$. On the other hand lower bound is given by 3. This follows from \cite[Corollary 2]{kucharz1991} (see also \cite[Question 4.1]{grimm2015}). In particular it is still not known if there exists an algebraic surface $X$ over $\R$ with $p(\R(X))=3$.

\section{Applications.}
In this section we will apply results obtained in Section 3 to uniformly rational surfaces

\subsection{Uniformly rational varieties}
Let $X \subset K^n$ be an irreducible algebraic set, such sets will be called  algebraic varieties. Take $Y \subset K^m$ to also be an algebraic variety.

\begin{df}
We say that $\varphi : X \rightarrow Y$ is a regular map if $\varphi=\left( \varphi_1, \dots, \varphi_m \right)$ where each $\varphi_i$ is a regular function on $X$. $\varphi$ is a biregular map if both $\varphi$ and $\varphi^{-1}$ are regular maps.
\end{df}

\begin{df}\label{rational}
We say that $X$ is a $K$-rational variety if there exists a Zariski open subset $U \subset X$ which is biregularly isomorphic with a Zariski open subset $V\subset K^m$ for some $m$.

We say that $X$ is a uniformly $K$-rational variety if every point of $X$ posses a Zariski open neighbourhood biregurarly isomorphic with a Zariski open subset of $K^m$.

\end{df}
If $X$ is a $K$-rational variety, then its field of fractions is isomorphic to $K(x,y)$. Any uniformly rational variety is nonsingular. Uniformly rational varieties have been studied for a long time and they are present in the literature under various names such as plain varieties \cite{bhsv2008}, regular \cite{gromov1989}, locally flattenable \cite{popov2020}.

It is a long standing open problem posed by Gromov \cite{gromov1989} whether a nonsingular rational variety is uniformly rational. Some instances of this problem are known \cite{bb2014}.
If $K$ is an arbitrary field then a blowup along a smooth center of a uniformly rational variety is uniformly rational \cite{bhsv2008}.
If $K=\R$ or $\mathbb{C}$, then it is known that any nonsingular $K$-rational curve is uniformly $K$-rational and any nonsingular $K$-rational surface is uniformly $K$-rational (cf. \cite{bk2025} for the real case and references therein). This problem is still open in full generality in dimension at least three.

Let $X_K \subset K^n$ be and algebraic set and let $\bar{K}$ be an algebraic closure of $K$. Denote by $X_{\bar{K}} \subset \bar{K}^n$ the algebraic set considered as an algebraic set over $\bar{K}$ defined by the same polynomials as $X_K$.

\begin{remark}
One has to be careful when dealing with various notions of rationality. Such ideas are also considered from the scheme theoretic point of view. Usually, an algebraic variety is called rational if $X_{\bar{K}}$ is rational - sometimes this property is called geometrical rationality. This differs substantially from Definition \ref{rational}. In particular, it may happen that an algebraic variety is not $K$-rational but it is geometrically rational (see \cite[Example 1.3]{kollar1995}). Definition \ref{rational} corresponds to the $K$-rationality \cite[Definition 2.1]{kollar1995} or $K$-birational triviality \cite[Definition 12.4]{manin1986}. 
\end{remark}
 Of course, if $X$ is $K$-rational then it is geometrically rational. 
 It is not clear what is the exact relation between geometrical rationality and $K$-rationality (cf. Section \ref{sec:open problems}). Here, we are able to show the following.

\begin{thm}
Let $X_K$ be a $K$-rational variety of dimension $n$ such that $X_{\bar{K}}$ is uniformly $\bar{K}$--rational variety. Then there exists a field $L$ such that $K 
\subset L \subset \bar{K}$ and $[L:K] <\infty$ such that $X_L$ is uniformly $L$--rational.

\end{thm}
\begin{proof}
Since $X_{\bar{K}}$ is uniformly rational, we may take a Zariski open covering $U_i$ for $i=1,2,\dots, m$ such that each $\varphi_i: U_i \rightarrow V_i \subset \bar{K}^n$ is a biregular isomorphism. As $\bar{K}$ is an algebraic extension of $K$, each of the aforementioned isomorphisms are given by a tuple of rational functions. Each rational function is a quotient of two polynomial functions from $\bar{K}[X_{\bar{K}}]=K[X]\otimes_K \bar{K}$, and there are only finitely  many elements from $\bar{K}\setminus K$ involved. Hence, for $L$ we may take take smallest field extension of $K$ such that all of the $\varphi_i$ and $\varphi_j^{-1}$ are defined over $L$ for $i,j=1,2,\dots, m$.

\end{proof}

Let us now move to the theory of sums of squares on uniformly rational varieties.
\begin{df}
We say that a uniformly $K$-rational variety of dimension $l$ is of type $m$, if $m$ is the smallest positive integer such that it posses a covering by Zariski open subset isomorphic to Zariski open subset of $K^l$.

\end{df}
$X$ is of type 1 iff $X \subset K^n$. If $K=\R$ then the unit sphere is of type 2 and the real projective space $\mathbb{P}^2(\R)$ is of type at most 3.
For a positive integer $k$ we define $\lceil k \rceil$ to be the smallest positive integer greater than or equal to $k$.  Denote by $f_d(m,n)=\left\lceil\frac{nm}{d} \right\rceil \left\lceil\frac{m}{d}\right\rceil d$ for $d=2,4,8$. By $K(x,y)$ we denote the field of rational functions in two variables. The second main result of this paper is the following.

\begin{thm}
Let $K$ be a formally real field. Assume that $p(K(x,y))=n < \infty$. If $X$ is a uniformly $K$-rational surface of type $m$ then the following inequality holds
$$p(\OO(X)) \leq \min\{f_2(m,n),f_4(m,n),f_8(m,n) \}. $$
In particular for large values of $m,n$ we see that $p(\OO(X))$ is at most of order $\frac{m^2n}{8}$.
\end{thm}
\begin{proof}
Fix a covering $U_i \subset X$ for $i=1,2,\dots, m$ such that each $U_i$ is isomorphic to $V_i \subset K^2$. Let $f \in \OO(X)$ be a sum of squares. Choose functions $h_i \in K[X]$ such that $\mathcal{Z}(h_i)=X\setminus U_i$ for $i=1,2,\dots, m$. Consider now $f_i:=f|_{U_i}$. By Corollary \ref{cor:p(O(X))=p(O(U))} we obtain a presentation 
$$f_i=\frac{\sum_{j=1}^n f^2_{ij}}{g_i^2}=\frac{\sum_{j=1}^n f^2_{ij}h_i^2}{g_i^2h_i^2}$$
where $\mathcal{Z}(g_i) \subset X\setminus U_i$ and $f_{ij},g_i \in K[X]$. As $\{U_i\}$ is a covering of $X$, we obtain in this manner a nowhere vanishing polynomial function $\sum_{i=1}^mg_i^2h_i^2$. By a standard trick used in real algebraic geometry, we may now add all of the numerators and denominators of all presentations to obtain a global one i.e.
$$f=\frac{\sum_{i=1}^m \sum_{j=1}^n f^2_{ij}h_i^2}{\sum_{i=1}^mg_i^2h_i^2}=\frac{\left(\sum_{i=1}^m \sum_{j=1}^n f^2_{ij}h_i^2\right)\left(\sum_{i=1}^mg_i^2h_i^2\right)}{\left(\sum_{i=1}^mg_i^2h_i^2\right)^2}.$$
This gives an upper bound on the length of $f$ as $m^2n$, in particular $p(\OO(X)) \leq m^2n$.
For the second part of the proof, recall that there are 2-, 4- and 8-square identity, in the sense of that
$(\sum_{i=1}^n x_i^2)(\sum_{i=1}^n y_i^2)=\sum_{i=1}^n z_i^2$
holds where $z_i$ is bilinear in $x_j$'s and $y_k$'s and $n=1,2,4,8$. Of course, there are such identities for every $2^k$ but the $z_i$'s are rational functions rather than bilinear in $x_j$'s and $y_k$'s.

Back to the proof, the numerator of $f$ can be written as a product of $nm$ and $m$ squares. We may now divide $nm$ and $m$ into the sets of $8$ squares each. In case $nm$ and $m$ are not divisible by $8$, we add 1 additional set, which we can fill with zeros if needed. By a repeatedly using the $8$ square identity we reduce the number of summands from $m^2n$ to $ \left\lceil\frac{nm}{8} \right\rceil \left\lceil\frac{m}{8}\right\rceil8 $ which is precisely $f_8(m,n)$. The same reasoning can be applied with the $4$-square and $2$-square identities which in total gives $p(\OO(X)) \leq \min\{f_2(m,n),f_4(m,n),f_8(m,n) \}$, finishing the proof.

\end{proof}

 If $m,n$ are large, then clearly $f_8(m,n)$ is the best bound we can get, however if the values of $m,n$ are small then it may happen that using the 4-square formula instead of $8$-square formula yield a smaller upper bound. For the sake of next section, let us consider the case $m=3$, i.e. uniformly rational varieties of type $3.$ Denote by $f_d(n)=f_d(3,n)$ for $d=2,4,8$. The precise values of these functions depend on the value of $n$ modulo 8. Below, we present a comparison of these functions in terms of $n$. Let $n=8k+l$ where $l=0,1,\dots, 7$.

\begin{center}
\begin{tabular}{ | m{1em} | m{1.5cm}| m{1,5cm} | m{1,5cm} | } 
  \hline
  $l$ & $f_2(n)$ & $f_4(n)$ & $f_8(n)$ \\ 
  \hline
 0 1 2 3 4 5 6 7  & 48k 48k+8 48k+12 48k+16 48k+24 48k+32 48k+36 48k+44 & 24k 24k+4 24k+8 24k+12 24k+12 24k+16 24k+20 24k+24 & 24k 24k+8 24k+8 24k+16 24k+16 24k+16 24k+24 24k+24 \\ 
 \hline
\end{tabular}
\end{center}

If $X$ is a uniformly rational variety over a field $K$, we have show the existence of a bound on $p(\OO(X))$ in terms of its type and $p(K(x,y))$. However in some cases much better result can be obtained.

\subsection{Projective space, sphere and torus.}
Let $K$ be a formally real field. In this case, it is known \cite[Theorem 3.4.4]{bcr1998} that the projective space is an affine algebraic set, hence we will not speak about projective algebraic sets, additionally in order to rephrase results about projective varieties we have to make some additional assumptions.

In this section we will consider blowups of three particular surfaces, namely the projective space $\mathbb{P}^2(K)$, the unit sphere given as a quadric $\{x^2+y^2+z^2=w^2 \}\subset \mathbb{P}^3(K)$ and $\mathbb{P}^1(K) \times \mathbb{P}^1(K)$ which is biregular to the torus $\mathbb{S}_K^1\times \mathbb{S}_K^1$. The isomorphism between $\mathbb{P}^1(K)$ and $\mathbb{S}_K^1$ is given by the stereographic projection. Here $\mathbb{S}_K^1$ is the unit circle in the plane $K^2$. Let us refer to the three aforementioned varieties, as \textit{basic rational surfaces}. By a multiblowup we will mean a finite composition of blowups in points. Let us now state a slight reformulation of a very useful fact:
\begin{prop}\cite[Lemma 1.6]{cssv2023}
Let $X$ be a nonsingular $K$-rational surface such that there exist three open subsets $U_0,U_1,U_2 \subset X$ with 
\begin{enumerate}
\item $U_0\cup U_1\cup U_2 = X$,
\item for all $i=0,1,2$ $U_i$ is isomorphic to $K^2$,
\item $U_i \setminus (U_j\cup U_k)$ is finite for any choice of $i,j,k=0,1,2$.
\end{enumerate}
Consider a finite set $A_1 \subset U_0 \cap U_1 \cap U_2$. Let $P \in  (U_0 \cap U_1 \cap U_2) \setminus A_1$ be a point and consider $\pi : Y \rightarrow X$ to be the blowup of $X$ at $P$. Consider also a finite set $A_2$ in the exceptional divisor $E=\pi^{-1}(P) \subset Y$. Then, there exist three open sets $U_0',U_1',U_2'\subset Y$ such that:

\begin{enumerate}
\item $U_0'\cup U_1'\cup U_2'=Y$
\item for all $i=0,1,2$ $U_i$ is isomorphic to $K^2$,
\item both $A_2$ and the proper transform of $A_1$ are contained in $ U_0'\cap U_1'\cap U_2'$.
\item $U_i' \setminus (U_j'\cup U_k')$ is finite for any choice of $i,j,k=0,1,2$.

\end{enumerate}

\end{prop}
The above is proved over the field of complex number \cite{cssv2023}, yet the proof works over an infinite field $K$. Basically, apart from the standard properties of a blowup, one uses infiniteness of $K$ and the property that $U_i\setminus ( U_j\cup U_k)$ is finite for every choice of $i,j,k=0,1,2$. Having this result at hand we may now state the following.

\begin{thm}
Let $Y$ be a multiblowup of a basic rational surface $X$. Then the type of $Y$ is at most 3.
\begin{proof}
One can repeat proof of \cite[Theorem 1.1]{cssv2023} verbatim.
\end{proof}
\end{thm}

This easily obtained the generalization of \cite[Theorem 1.1]{cssv2023} can be now combined with a previous section. 

\begin{thm}
Let $p(K(x,y))=n<\infty$ and let $Y$ be a multiblowup of a basic rational variety, or a Zariski open subset of such. Then $p(\OO(Y))\leq \min \{f_2(n),f_4(n),f_8(n) \}$
\end{thm}

\begin{remark}
Having the above results at hand, one question remains widely open, namely: for which fields $K$ we have $p(K(x,y)) < \infty$? Even more, it is still not known whether $p(K)<\infty $ implies $p(K(x)) < \infty$. Some partial results are known. If $K$ is a real closed field then $p(K(x,y)) \leq 4$  \cite[Corollary 3.4]{pfister1995}, if $K$ is a (hereditarily) pythagorean field then $p(K(x,y)) \leq 8$ \cite{bdgmfz2023}. If $K$ is a number field then $p(K(x,y)) \leq 8$ \cite[Introduction]{bdgmfz2023}.

\end{remark}

In particular if $K=\R$ we obtain two important corollaries from the above.

\begin{thm}
If $Y$ is a nonsingular compact $\R$-rational algebraic surface, then $p(\OO(Y))\leq 12$.

\begin{proof}

It is known that any nonsingular compact $\R$-rational algebraic surface is a multiblowup of $\mathbb{P}^2(\R)$, a sphere, or $\mathbb{P}(\R)\times \mathbb{P}(\R)$ (cf. \cite[Theorem 3.1]{bh2007}). Hence the result.

\end{proof}
\end{thm}

Of course, if $ U\subset Y$ is a Zariski open subset of such a surface, then $p(\OO(U))\leq 12$. If we further specialize to the case where $Y$ is a multiblowup of $\R^2$ then

\begin{cor}
If $Y$ is a multiblowup of $\R^2$ then $p(\OO(Y)) \leq 12$.
\end{cor}

If $K=\R$ then the real projective space is actually a blow up of a sphere in 1 point. 
\begin{remark}
In the above corollary, the more important is the uniform bound, rather that the bound is 12. Recently there were several papers concerning the various problems for $k$-regulous functions \cite{banecki2025, bk2023}. The author, together with Banecki \cite{bk2023}, have shown that if $X\subset \R^n$ is a closed irreducible $0$-regulous set, then $p(\RR^0(X))\leq 2^n$ where $n$ is the dimension of $X$ and $\RR^0(X)$ is the ring of $0$-regulous functions on $X$. Also, every $0$-regulous function which is nonnegative is a sum of squares. Such statement is not true for $k>0$. Not much is known about the Pythagoras numbers of $k$-regulous functions on $\R^2$ for $k>0$. However, the $k$-regulous functions on $\R^n$ are precisely those continuous functions of class $\mathcal{C}^k$ which become regular after a suitable multiblowup. By the above corollary, having a uniform bound on the Pythagoras number of the ring of regular functions on blowups of $\R^2$ it may be possible to somehow determine if the Pythagoras number of the ring of $k$-regulous functions on $\R^2$ are finite or not.
\end{remark}

\begin{remark}
If $X$ is a nonsingular rational surface over $\R$ then we showed the finiteness of the Pythagoras number of $\OO(X)$. This implies finiteness of the unital Pythagoras number cf. \cite[Theorem 5.17]{bp1996}. This on the other hand implies finiteness of the higher unital Pythagoras numbers. This means, that there is an explicitly computable uniform upper bound $N_n$ such that any unit in $\OO(X)$ which is a sum of $2n$th powers is a sum of at most $N_n$ $2n$th powers. This can be seen as a generalization of the Theorem of Mah\'e \cite[Theorem 7.3]{mahe1990} to higher powers. As mentioned in the introduction, this is far from computing the higher Pythagoras numbers, as we have no control on the lengths of non--units. In particular, the  $2n$-length of polynomials constructed in \cite{kv2023} does not tend to infinity in the ring $\OO(\R^2)$.

\end{remark}

\section{Further remarks}
It this section we discuss the possible generalization of the Theorem \ref{main} to a larger class of surfaces as well as to higher powers.
\subsection{Lack of Theorem of Cassels for coordinate rings of curves.}
It would be valuable to generalize Theorem \ref{main} to the rings of regular functions of products of two irreducible nonsingular curves. This would easily follow from Theorem of Cassels for coordinate rings of curves. However, such a generalization is false as the following example shows. Let $K=\mathbb{R}$ and $k \in \mathbb{R}[x]$ be a nonlinear polynomial of odd degree in a single variable $x$. Denote by $X=\{y^2-k(x)=0 \} \subset \R^2$ a  planar curve.

\begin{prop}
With the above notation:
\begin{itemize}
\item[a)] The Pythagoras number of $\R[X]$ is 3 or 4.
\item[b)] The Pythagoras number of the field $\R(X)$ is 2.

\end{itemize}

\begin{proof}
Part $a)$ is precisely \cite[Proposition 3.8]{cldr1982}. For the part $b)$, observe that $\R(X)$ has transcendence degree 1 over $\R$, hence by the result of Pfister \cite[Corollary 3.4]{pfister1995} we have $p(\R(X))\leq 2$. Straightforward calculations shows that $1+\bar{x}^2$ is not a square in $\R(X)$, where $\bar{x}$ is the image of $x$ in $\mathbb{R}[X]$.

\end{proof}
\end{prop}

With the above at hand, let $\varphi = x_1^2+x_2^2$ be a quadratic form. There exists an element $f\in \R[X]$ which is a sum of at least three squares but not of two squares. On the other hand $\varphi$ represents this element in the fraction field. Hence, Cassels Theorem cannot be generalized to coordinate rings of curves, not even nonsingular ones.

\subsection{Lack of Theorem of Cassels for higher powers.}
Let $n>1$ be a positive integer. Theorem of Cassels also does not hold for sums of higher powers instead of squares. This follows from the fact that there are polynomials in a single variable $x$ which are sums of $2n$-th powers of rational functions, but no sums of $2n$-th powers of polynomials \cite{prestel1984}. Let us discuss the case of fourth powers, as examples for higher powers can be produced similarly. 

Consider the family of polynomials $f_m=x^4+mx^2+1$ \cite[Exercise 7.4.4]{pd2001}. An  application of Cauchy-Schwartz-Bunyakovsky inequality shows that $f_m$ is a sum of fourth powers of polynomials only for $m \in [0,6]$. On the other hand, by the Theorem of Becker (see \cite{becker1982, clpr1996}) $f_m$ is a sum of fourth powers of rational functions for $m \in [0, +\infty)$.

\section{Open problems}\label{sec:open problems}
We finish with a section devoted to some open problems.

In  Section 2 we have discussed the notion of weakly factorial varieties. We gave a topological classification of those over the real numbers. No similar characterization is known over an arbitrarily field of characteristic zero. Hence, we propose the following problem:

\begin{prob}
Find a characterization of weakly factorial varieties over an arbitrary field of characteristic zero.

\end{prob}
In the Section 3 we obtained a strengthening of a result of Scheiderer. However, we believe that a much stronger statement is true.
\begin{conj}
Let $R$ be a regular local ring of dimension 2 and $K$ be its field of fractions. Then $p(R)=p(K)$.

\end{conj}
Of course, $p(R)\geq p(K)$, hence this conjecture is only interesting if $p(K)$ is finite.

The main ingredient in the proof of Theorem \ref{main} was Theorem of Cassels. As we have seen, such theorem does not hold for higher powers. Hence, the proof of Theorem \ref{main} does not generalize to higher powers. The relation between $p_{2d}(\OO(X))$ and $p_{2d}(K(X))$ is not yet known. Here, $p_{2d}$ is the $2d$-Pythagoras number (cf. \cite{kv2023})

\begin{prob}
Study sums of higher powers and higher Pythagoras numbers in the rings of regular functions on algebraic surfaces over a formally real field $K$.

\end{prob}

Next question relates the notions of $K$ and $\bar{K}$ uniform rationality.
\begin{conj}
Let $X_K \subset K^n$ be a nonsingular $K$-rational algebraic variety. Assume that $ X_{\bar{K}} \subset \bar{K}^n$ is uniformly $\bar{K}$-rational. Is $X_K$ uniformly $K$-rational?
 \end{conj}
We have shown in Theorem 4.4 that under the above assumptions, there exists a field extension $L$ of $K$ of a finite degree
such that $X_L$ is uniformly $L$-rational. The general conjecture seems feasible.
If the answer to the above conjecture is negative, then it may be interesting to study field extensions of minimal degree with respect to the above property. By Theorem 4.4 such fields exist. 
The assumption of $X$ being rational is necessary, as $\bar{K}$-rationality does not imply rationality over $K$.

We finish this section with a conjecture concerning (rational) surfaces.
\begin{conj}
Let $X$ be a nonsingular ($K$-rational) surface over a formally field $K$. Then $p(\OO(X))=p(K(X))$.

\end{conj}
Result as above would imply similar result for large class of 2-dimensional local rings.

\section*{Acknowledgment}
The author is very grateful to Bartłomiej Bychawski, Andrzej Czarnecki and Piotr Miska for their comments, remarks and many fruitful discussions. Also, to Juliusz Banecki for pointing out the application of Theorem \ref{main} to the uniformly rational varieties in a previous version of this paper.

\bibliographystyle{plain}
\bibliography{refs.bib}

\begin{small}

\vspace{5pt}

\noindent
Tomasz Kowalczyk

\noindent
Institute of Mathematics

\noindent
Faculty of Mathematics and Computer Science

\noindent
Jagiellonian University

\noindent
ul. Łojasiewicza 6, 30-348 Kraków, Poland

\noindent
e-mail: tomek.kowalczyk@uj.edu.pl

\end{small}

\end{document}